\documentclass[a4paper,12pt,oneside,reqno]{amsart}
\usepackage[headinclude,DIV13]{typearea}
\areaset{15.5cm}{25cm}
\usepackage{amsfonts,amssymb,amsmath,amsthm}

\newtheorem{theorem}{Theorem}[section]
\newtheorem{lemma}[theorem]{Lemma}
\newtheorem{proposition}[theorem]{Proposition}
\newtheorem{corollary}[theorem]{Corollary}

\theoremstyle{definition}

\theoremstyle{remark}
\newtheorem*{remark}{Remark}

\numberwithin{equation}{section}

\def\psum{\mathop{\sum\nolimits'}}
\def\d{\mathrm{d}}

\def\bbone{{\mathchoice {\rm 1\mskip-4mu l} {\rm 1\mskip-4mu l}
          {\rm 1\mskip-4.5mu l} {\rm 1\mskip-5mu l}}}

\def\<{\langle}
\def\>{\rangle}

\begin{document}
\title{Hitting time of large subsets of the hypercube}

\author[J. \v Cern\'y]{Ji\v r\'\i~\v Cern\'y}
\address{Ji\v r\'\i~\v Cern\'y\\
  \'Ecole Polytechnique F\'ed\'erale de Lausanne\\
  1015 Lausanne\\
  Switzerland
  }
\email{jiri.cerny@epfl.ch}

\author[V. Gayrard]{V\'eronique Gayrard}
\address{V\'eronique Gayrard\\Laboratoire d'Analyse, Topologie, Probabilit\'es\\
CMI, 39 rue Joliot-Curie\\
13453 Marseille Cedex}
\email{gayrard@latp.univ-mrs.fr, veronique@gayrard.net}


\date{\today}
\begin{abstract}
  We study the simple random walk on the $n$-dimensional hypercube, in
  particular its hitting times of large (possibly random) sets. We give
  simple conditions on these sets ensuring that the properly-rescaled
  hitting time is asymptotically exponentially distributed, uniformly in
  the starting position of the walk. These conditions are then verified for
  percolation clouds with densities that are much smaller than
  $(n \log n)^{-1}$. A main motivation behind this paper is the study of
  the so-called aging phenomenon in the Random Energy Model (REM),
  the simplest model of a mean-field spin glass. Our results allow
  us to prove aging in the REM for all temperatures,
  thereby extending earlier results to their optimal temperature domain.

\end{abstract}

\maketitle

%%%% START OF ht1.tex %%%%
\section{Introduction}
\label{s:introduction}
Let $\mathcal V_n$ be the $n$-dimensional hypercube,
$\mathcal V_n=\{0,1\}^n$. We equip $\mathcal V_n$ with
the metric
\begin{equation}
  d(x,y)=\sum_{i=1}^n \bbone\{x(i)\neq y(i)\},
\end{equation}
where $x(i)$ are the coordinates of $x\in \mathcal V_n$. Let further $Y_n$ be
the simple random walk on $\mathcal V_n$. That is,  $Y_n$ is the
discrete-time Markov chain with state space $\mathcal V_n$ whose
transition probabilities are given by
$\mathbb P[Y_n(k+1)=y|Y_n(k)=x]=n^{-1}$ if $d(x,y)=1$, and $0$ otherwise. We
write $\mathbb P_x$ for the distribution of $Y_n$ conditioned on $Y_n(0)=x$.
For $A\subset \mathcal V_n$ we define the hitting time of $A$ by
\begin{equation}
  H_n(A)=\min\{k\ge 0: Y_n(k)\in A\}.
\end{equation}

We are interested in the distribution of the hitting time of large random
subsets of the hypercube. Specifically, let $\rho \in [0,1] $. We say that
the set $A\subset \mathcal V_n$ is a \textit{percolation cloud} on
$\mathcal V_n$ with density $\rho $ if  each site $x\in \mathcal V_n$ is
in $A$ with probability $\rho $ independently of all others.

Our main aim is to prove the following theorem.
\begin{theorem}
  \label{t:perc}
  Let $\bar m(n)$ be such that
  \begin{equation}
    \label{e:condm}
    n \log n \ll \bar m(n) \ll 2^{n}(\log n)^{-1}
  \end{equation}
  and let $A_n$ be a sequence of percolation clouds on $\mathcal V_n$ with
  densities $\bar m(n)^{-1}$ defined on a common probability space
  $(\Omega,F, P)$.
  Then, for all $a>0$,
  \begin{equation}
    \label{e:thmjedna}
    \lim_{n\to\infty}
    \max_{x\in \mathcal V_n}
    \Big |
    \mathbb P_x \big[H_n(A_n\setminus x)\ge a \bar m(n)\big]-e^{-a}
    \Big |
    =0\,,\quad P\text{-a.s.}
  \end{equation}
  In words, the distribution of the normalised hitting time of
  $A_n\setminus x$ converges to the exponential distribution uniformly in
  the starting position $x$.
\end{theorem}

In Theorem~\ref{t:sampling} below we state a similar result for
another important class of random sets, namely sets $A_n$ that are
sampled from $\mathcal V_n$ without replacement (i.e. each subset
of $|A_n|$ elements of $\mathcal V_n$ is equally likely).

\begin{theorem}
  \label{t:sampling}
  Let $M_n$ be a sequence of integers such that, setting $m(n)={2^n}/{M_n}$,
  \begin{equation}
    \label{e:condm'}
     n \log n \ll m(n) \leq 2^n
  \end{equation}
  Let $A_n$ be subsets of $M_n$ elements sampled from $\mathcal V_n$ without replacement,
  and defined on a common probability space
  $(\Omega',F', P')$.
  Then, for all $a>0$,
  \begin{equation}
    \lim_{n\to\infty}
    \max_{x\in \mathcal V_n}
    \Big |
    \mathbb P_x \big[H_n(A_n\setminus x)\ge a m(n)\big]-e^{-a}
    \Big |
    =0\,,\quad P'\text{-a.s.}
  \end{equation}
\end{theorem}

Estimates on the distribution of the hitting time of various subsets of
the hypercube have a long history. They can be traced back to the early
literature on first passage times for Markov chains (see \cite{Kem61} and
  the references therein) where these questions are reformulated in terms
of the (one-dimensional) Ehrenfest urn scheme. These results provide very
sharp estimates on the asymptotic distribution of the hitting time of a
single point. More, recently Matthews \cite[p.118]{Mat89} gave finite-$n$
estimates for the Laplace transform of the hitting time of sets containing
one or two points.  These estimates are key ingredients that enter his
description of the covering time of the hypercube and related questions.

The hitting times of much more general sets (possibly random, and whose
  size is possibly increasing with $n$) were studied very recently in
\cite{BC06} and \cite{BG06}. In \cite{BG06}, Ben Arous and Gayrard give
precise conditions for the hitting time of subsets of $\mathcal V_n$ to be
asymptotically exponentially distributed for a class of subsets  for which
the so-called lumping construction can be applied.  This construction was
fist used in this context by \cite{BBG03}. It can be understood as a $d$-dimensional
extension of the Ehrenfest urn scheme where the random walk $Y_n$ on the
hypercube is replaced by a walk defined on a $d$-dimensional state space
of smaller cardinality, and which evolves in a convex potential that is
very steep near its boundary. Such a chain is then studied using the tools
developed in \cite{BEGK1},\cite{BEGK2} for the study of metastability in
reversible Markov chains on discrete state space. This method allows in
particular to show that Theorem~\ref{t:perc} is valid for
$m(n)\ge C 2^n/\log n$.

In \cite{BC06}, Ben Arous and \v Cern\'y obtain a result similar to
Theorem~\ref{t:perc}. Namely, they prove convergence to the
exponential distribution of the hitting time of percolation clouds
for densities of order $2^{-c n}$ with $c\in (3/4,1)$ (see
  Lemma~3.7 of \cite{BC06}). The method is based on the formula discovered
by Matthews in his study of covering problems (see Theorem 1.3 in
  \cite{Mat88}) and on (improved) estimates from \cite{Mat89}.

  The results of the present paper are the first that allow to treat sets
 of very large size, namely sets as large as $o(2^n (n \log n)^{-1})$.

\smallskip

We will show in Section~\ref{s:equiv} that
both Theorem~\ref{t:perc} and Theorem~\ref{t:sampling} are
consequences of the following more general theorem. To state it we define
the function $\xi_n(k)$ by
\begin{equation}
  \label{e:defxi}
  \xi_n(k)=2^{-n } \frac n 2 \binom n k^{-1}
  \sum_{j=1}^{n-k}\binom n {k+j} \frac 1 j,
  \qquad k=0,\dots ,n.
\end{equation}
The role of this relatively complicated function will become evident in
Lemma~\ref{l:twopoint}.

\begin{theorem}
  \label{t:gen}
  Let $n \log n \ll m(n) \leq 2^n$ and
  $B_n\subset \mathcal V_n$ be such that
  \begin{equation}
    \label{e:Bsize}
    0<|B_n|= 2^n m(n)^{-1}(1+o(1)).
  \end{equation}
  Define
  \begin{equation}
    \begin{split}
      \label{e:volumes}
      v_n(k)&=\max_{x\in \mathcal V_n} \big|\{y\in B_n:d(x,y)=k\}\big|,\\
      V_n(k)&=\max_{x\in \mathcal V_n} \big|\{y\in B_n:d(x,y)\le k\}\big|.
    \end{split}
  \end{equation}
  If there exists a function $g(n)$, such that $g(n)\le n/2$,
  \begin{equation}
    \label{e:gg}
    \xi_n(g(n))\ll 2^{-n}m(n),
  \end{equation}
  and
  \begin{equation}
    \label{e:vcond}
    \lim_{n\to \infty}\sum_{k=1}^{g(n)-1}v_n(k)\xi_n(k)=0,
    \qquad
    \text{and}
    \qquad
    V_n(g(n)-1)\ll |B_n|,
  \end{equation}
  then
  \begin{equation}
    \label{e:lim}
    \lim_{n\to\infty}
    \max_{x\in \mathcal V_n}
    \Big |
    \mathbb P_x \big[H_n(B_n\setminus x)\ge a m(n)\big]-e^{-a}
    \Big |
    =0.
  \end{equation}
\end{theorem}

\begin{remark}
  1. While the upper bound on $m(n)$ in the range of validity of the
  preceding theorem is trivially optimal, we cannot claim the same about
  the lower bound. Morally, the proof exploits the fact that the simple
  random walk on the hypercube equilibrates after $O(n \log n)$ steps,
  even if an application of this fact is not easy to see in our actual
  proof. It is intuitively obvious that after equilibration the hitting
  time should be exponentially distributed. The claim of the theorem might
  stay valid even for $m(n)=O(n\log n)$, however for different reasons.
  Several serious technical complications appears in our proof for such
  $m(n)$.  We were not motivated to proceed further since it is not needed
  for the applications we have in mind.

  2. The results of  \cite{BG06} (see Theorem 1.7 of \cite{BG06} and the
    simpler Corollary 1.8) imply that Theorem~\ref{t:perc} holds true for
  sets $B_n$ whose size is asymptotically constant that is, when
  \eqref{e:Bsize} holds for $m(n)={2^n}/{M}$ where $M$ is a fixed integer.

  3. We will show soon that $\xi_n(1)=O(n^{-1})$. Therefore, the
  conditions appearing in \eqref{e:vcond} are void when \eqref{e:Bsize}
  holds with $m(n)$ satisfying $2^{-n}m(n)\gg n^{-1}$ (and, in
    particular, when $|B_n|$ is constant).  As an immediate consequence of
  Theorem~\ref{t:gen} we then obtain:
\end{remark}

\begin{corollary}
  \label{c:perc}
  Theorem~\ref{t:perc} remains valid on $2^{n}(\log n)^{-1}\le \bar m(n)\leq
  2^n$ if the term $a\bar m(n)$ in \eqref{e:thmjedna} is replaced by $2^n |A_n|^{-1}$.
\end{corollary}

\begin{proof}
  It is sufficient to take $m(n)$ in Theorem~\ref{t:gen}
  to be $m(n)=2^n/|A_n|$, where $A_n$ is the percolation cloud with
  intensity $\bar m(n)$.
\end{proof}

Our method of proof relies on a sharp estimate on the distribution
of hitting time of a single point when time is measured on the
scale $m(n)$, namely on the quantity
\begin{equation}
  \mathbb P_{z_k}[H_n(\boldsymbol 0)< a m(n)]
\end{equation}
where $\boldsymbol 0$ is the vertex of $\mathcal V_n$ whose coordinates
are all $0$ and $z_k$ is any vertex at distance $k$ from it; this estimate
is itself derived from an estimate on the Laplace transform of the hitting
time $H_n(\boldsymbol 0)$ (see respectively Lemma~\ref{l:twopoint} and
  Lemma~\ref{l:laplace}). The domain of validity of the latter determines
our bounds on $m(n)$.  On this domain the Laplace transform is well
approximated by the sum of two terms:  the expected contribution of an
exponential distribution, and the mysterious looking function $\xi_n(k)$
which is, essentially, the probability that $Y_n$ started from $z_k$ hits
$\boldsymbol 0$ in the first $n\log n$ steps.

The rest of the paper is organised as follows. In the next section we
describe our main motivation and give some important consequences of our
results for aging in the Random Energy model.
In Section~\ref{s:prop} we give the proof of Theorem~\ref{t:gen}.
Finally, Theorems~\ref{t:perc} and~\ref{t:sampling} are proved in
Section~\ref{s:equiv}.
%%%% END OF ht1.tex %%%%

%%%% START OF ht4.tex %%%%

\section{Aging in the Random Energy Model}
\label{s:aging}

The main motivation behind this paper originates in the study of the
Random Hopping Time (RHT) dynamics of the Random Energy Model (REM), which
is often called the simplest model of a spin glass. Let us describe this problem
briefly (for a recent review  see \cite{BC06c}). In the REM an energy $E_x$
is attached to every site $x\in \mathcal V_n$. The $E_x$'s are chosen to
be i.i.d.~with standard normal distribution. Given a collection
$\boldsymbol E=\{E_x:x\in \mathcal V_n\}$ and a parameter $\beta>0 $
(representing the inverse of the temperature), the RHT dynamics in the REM
is defined as a continuous-time Markov process $X_n(\,\cdot\,)$ whose
transition rates are given by
$w^{\boldsymbol E}_n(x,y)=\exp(-\beta \sqrt n E_x)$ if $d(x,y)=1$ and zero
otherwise.

The main goal of the study of the processes $X_n$ was to prove aging. In
this context it usually means  showing that
the  two-point function
\begin{equation}
  R_n(t_w,t_w+t;\boldsymbol E):=P[X_n(t_w+t)=X_n(t_w)|\boldsymbol E],
\end{equation}
has a non-trivial limit as $t_w$, $ t=\theta t_w$ and $n$ tend simultaneously
to infinity.

The first proof of aging in the REM for $\beta > \beta_c=\sqrt {2 \log 2}$
and times $t_w\sim \exp (\beta \beta_c n)$ was given in
\cite{BBG03,BBG03b}, based on the analysis of the metastable behaviour of
$X_n$ and renewal theory. In \cite{BC06} another proof
of aging, based on the arc-sine law for stable subordinators, was given
for temperatures satisfying
\begin{equation}
  \label{e:bc}
  \sqrt {3/4} < \alpha \beta / \beta_c < 1.
\end{equation}
and (shorter) time-scales
$t_w(n)=( \alpha \beta \sqrt{2 \pi n})^{-1/\alpha }\exp(\alpha \beta^2 n)$,
where $\alpha \in (0,1)$ is a free parameter. We will now explain how our
methods allows to improve the lower bound $\sqrt{3/4}$ in \eqref{e:bc} to its
optimal value, that is to $0$. To this end we need to describe briefly the
background of the techniques of \cite{BC06}.

The transition rates of the process $X$ do not depend on the energy of the
target vertex $y$. Therefore, the process $X_n$ is  a time change of the simple
random walk $Y_n$, and  can be written as
$X_n(t)=Y_n(S^{-1}_n(t))$, where
\begin{equation}
  S_n(k)=\sum_{i=0}^{k-1} e_i \exp(\beta \sqrt n E_{Y_n(i)}),
\end{equation}
$S^{-1}_n$ is the generalised right-continuous inverse of $S_n(k)$, and
$e_i$ is a sequence of mean-one i.i.d.~exponential random variables.
$S_n(k)$ is the time of the $k^{\text {th}}$ jump of $X_n$.

It was proved in \cite{BC06} that $S(n)$ behaves (for $n$ large) as an
$\alpha $-stable subordinator  in certain time and temperature regimes. To
be more precise, set $r(n)=\exp(\alpha^2 \beta^2 n /2)$. Then, if
\eqref{e:bc} holds then   $t(n)^{-1}S_n(\lfloor  r(n) \,\cdot\,\rfloor)$
converges in distribution to an $\alpha $-stable subordinator.

It is a known fact that the value of a stable subordinator at time $t$ can
be approximated by the finite sum of its largest jumps up to this time.
The same is true for $S_n$: the main contribution to $S_n(r(n))$ comes
from a finite number of visits to sites with
$\exp(\beta \sqrt n E_x)\asymp t(n)$. Such sites  form (due to the
  i.i.d.~property of the energies) a percolation cloud.  To understand
properties of $S_n$ it is therefore necessary to understand how the simple
random walk visits such clouds.

As we have already remarked, the methods of \cite{BC06} are sufficient to
show Theorem~\ref{t:perc} for densities much smaller than $2^{-3n/4}$.
The constant $3/4$ in the exponent entails the constant $3/4$ in \eqref{e:bc}. The
result of our Theorem~\ref{t:perc} thus allows to extend the domain of
validity of Theorem 3.1 in \cite{BC06}.  As we find this extension important
we state it here:

\begin{theorem}
  For $\alpha \in (0,1)$ let $t_w(n)$ and $r(n)$ be as above. If
  \begin{equation}
    \label{e:new}
    0 < \alpha \beta < \beta_c,
  \end{equation}
  then $t_w(n)^{-1}S_n(r(n)\,\cdot\,)$ converges in distribution to an
  $\alpha $-stable subordinator and the two-point function $R_n$ exhibits
  aging. Namely, for a.e.~realisation of
  $\boldsymbol E$ and for every $\theta\in (0,\infty)$
  \def\Asl{\mathop{\mathsf{Asl}}\nolimits}
  \begin{equation}
    \lim_{n\to\infty} R_n(t_w(n),(1+\theta ) t_w(n);\boldsymbol E)=\Asl_\alpha (1/1+\theta ),
  \end{equation}
  where $\Asl_\alpha (u)$ stands for the distribution function of the
  generalised arcsine law with parameter $\alpha $,
  $ \Asl_\alpha (z):= \pi^{-1}\sin \alpha \pi  \int_0^{z}u^{\alpha
      -1}(1-u)^{-\alpha }\,\d u$.
\end{theorem}

It is worth noting that for any $\beta<\infty $ there exists $\alpha \in (0,1)$
such that \eqref{e:new} is satisfied. This implies that aging occurs in the
RHT dynamics of the REM at all temperatures.

\smallskip

Theorem~\ref{t:perc} allows further to study the RHT dynamics in
another interesting time-temperature regime where denser percolation
clouds should be considered:  for
$t_w(n)=(\beta \sqrt{2 \pi n})^{-1}\exp(\beta^2 n)$, that is for
$\alpha =1$. It is argued in the physics literature \cite{BB02b}
that the two-point function $R_n$ exhibits some interesting ultrametric
behaviour in this case. The rigorous treatment of this problem was the
original motivation behind this paper and will be subject of
a forthcoming paper.

%%%% END OF ht4.tex %%%%

%%%% START OF ht2.tex %%%%

\section{Proof of Theorem~\ref{t:gen}}
\label{s:prop}

As mentioned earlier Theorem~\ref{t:gen} is already known for small sets,  namely when
$m(n)={2^n}/{M}$ for $M$ a fixed integer (see \cite{BG06}, Corollary 1.8). Although
our method of proof clearly allows us to cover this case as well,
its treatment is, in places, quite different from the case $n \log
n \ll m(n) \ll 2^n$. Thus, in order to keep this paper as concise as
possible we will prove Theorem~\ref{t:gen} in that
latter case only and assume from now on that $n \log n \ll
m(n) \ll 2^n$.

\begin{proof}[Proof of Theorem~\ref{t:gen} (for $n \log n \ll m(n) \ll 2^n$)]
  We have
  \begin{equation}
    \mathbb P_x[H_n(B_n\setminus x)\ge a m(n)]
    =
    1- \mathbb P_x\Big[\bigcup_{y\in B_n\setminus x}
      H_n(y)< a m(n)\Big].
  \end{equation}
  Therefore,
  it follows from the inclusion-exclusion principle that
  for all even $\ell\in \mathbb N$
  \begin{equation}
    \mathbb P_x[H_n(B_n\setminus x)\ge a m(n)]
    \le
    1+\sum_{i=1}^\ell \frac {(-1)^i} {i!}
    \psum_{y_1,\dots,y_i\in B_n\setminus x}
    \mathbb P_x\Big[\bigcap_{j=1}^i H_n(y_j)< a m(n)\Big],
  \end{equation}
  where $\psum$ denotes the sum over all mutually different $y$'s.
  Analogous expressions for~$\ell$ odd give lower bounds.

  The following proposition is the
  key step of the proof.
  \begin{proposition}
    \label{p:sum}
    Under the assumptions of Theorem~\ref{t:gen}
    and assuming that $n \log n \ll m(n) \ll 2^n$ we have,
    for all $i\in \mathbb N$,
    \begin{equation}
      \lim_{n\to\infty}
      \max_{x\in \mathcal V_n}
      \bigg|
      \psum_{y_1,\dots,y_i\in B_n\setminus x}
      \mathbb P_x\Big[\bigcap_{j=1}^i H_n(y_j)< a m(n)\Big]
      -a^i\bigg| =0.
    \end{equation}
  \end{proposition}
  Using Proposition~\ref{p:sum} the completion of the proof of Theorem~\ref{t:gen}
  under the assumption that $n \log n \ll m(n) \ll 2^n$ is immediate.
\end{proof}

The proof of Proposition~\ref{p:sum} relies on several technical
lemmas which we collect in the subsection below. The proof of
Proposition~\ref{p:sum} is then concluded in subsection
\ref{s:propproved}.

\subsection{Preparatory Lemmas}
\label{s:preplem}

Our first lemma collects the properties of the function $\xi_n(k)$
which will be needed later.

\begin{lemma}
  \label{l:xi}
  (i) For all $k\in \{1,\dots,n\}$,
  $ \xi_n(k) \le K  \binom n k^{-1} n^{1/2}\log n$, where $K$ is a constant
  independent of $n$ and $k$.

  (ii) For all  $k\le n/2$,  $ \xi_n(k)\ge    \frac 1 2 \binom n k^{-1}$.

  (iii) For any fixed $n$ the function $\xi_n(k)$ is decreasing in $k$.

  (iv) If $k=o(n)$, then $ \xi_n(k)=\binom n k^{-1}(1+o(1))$.
\end{lemma}
\begin{proof}
  (i) Recall that
  $\xi_n(k)= 2^{-n}\frac n 2 \binom n k^{-1}\sum_{j=1}^{n-k}\binom n {j+k} \frac 1 j$.
  From a standard moderate deviations argument it follows that
  \begin{equation}
    2^{-n} \frac n 2 \sum_{\substack{j=0\\|j-n/2|\ge n^{7/12}}}^n \binom n j
    \le c_1 \frac n 2 e^{-c_2 n^{1/6}}.
  \end{equation}
  Therefore, the contribution of $j$'s with $|k+j-n/2|\ge n^{7/12}$ is
  $o(\binom  n k^{-1})$. For the remaining $j$'s we use the approximation
  \begin{equation}
    \label{e:binomass}
    \binom n {n/2+i} = \sqrt{\frac 2 \pi } n^{-1/2}2^n
    e^{-2 i^2/n}(1+o(1)),
  \end{equation}
  which is valid uniformly for $i=o(n^{2/3})$. Setting
  $a=n^{-1/2}(k+j-n/2)$ and $b=n^{-1/2}(k-n/2)$, and thus $j=n^{1/2}(a-b)$,
  the  contribution of the $j$'s with $|k+j-n/2|\le n^{7/12}$ to
  $\xi_n(k)$ equals
  \begin{equation}
    \begin{split}
      \label{e:gausapprox}
      2^{-n}\frac n2&
      \binom  n k^{-1}
      \sum_{\substack{a\in [-n^{1/12},n^{1/12}]\cap (\mathbb Z/\sqrt n)\\
          a\ge b+n^{-1/2}}}
      \binom n {a \sqrt n + n/2} \frac 1 {(a-b)\sqrt n}
      \\&\le K'
       n^{1/2}
      \binom  n k^{-1}
       \int_{( b+ n^{-1/2})\vee-n^{1/12}}^{n^{1/12}}
      e^{-2x^2}\frac {\d x} {(x-b)}
      \le
      K
      \binom  n k^{-1}
      n^{1/2}\log n.
    \end{split}
  \end{equation}

  (ii) For $k<n/2$,
  $\xi_n(k)\ge 2^{-n}\frac n 2 \binom n k^{-1} \sum_{j=1}^{n/2}\binom n {k+j}\frac 2 n
  \ge \frac 1 2 \binom n k^{-1}$.

  (iii) The function $\xi_n$ can be rewritten as
  \begin{equation}
    \label{e:secondform}
    \xi_n(k)=
    2^{-n}\frac n2
    \sum_{j=1}^{n-k}
    \binom{n-k} j
    \binom{k+j} j^{-1}
    \frac 1 j .
  \end{equation}
  Here, the fact that $\xi_n(k)$ is decreasing is apparent.

  (iv)  Using again the moderate  deviations argument, for $k=o(n)$,
  \begin{equation}
    \xi_n(k)=2^{-n}\frac n 2\binom n k^{-1}
    \sum_{j=n/2-n^{7/12}-k}^{n/2+n^{7/12}-k}\binom n {j+k} \frac
    1 j(1+o(1))=\binom n k^{-1}(1+o(1)).
  \end{equation}
  This completes the proof of Lemma~\ref{l:xi}.
\end{proof}

We now prove that for $m(n)$ in the range considered in Theorem~\ref{t:gen}
there always exists a function $g(n)$ with $g(n)\le n/2$ satisfying \eqref{e:gg}.
We in fact prove a little more:
\begin{lemma}
  Let $m(n)$ be such that $n \log n \ll m(n)\ll 2^n$. Then there exist
  function $g$ such that
  \begin{equation}
    \label{e:g}
    n/2-g(n)\gg n^{1/2}
    \qquad\text{and}\qquad
    \xi_n(g(n))\ll 2^{-n}m(n).
  \end{equation}
\end{lemma}
\begin{proof}
  Take $m'(n)$ such that $n\log n \ll m'(n)\ll m(n)$ and define $g(n)$ by
  \begin{equation}
    g(n)=\min\{ k: \xi_n(k)\le 2^{-n}m'(n)\}.
  \end{equation}
  Such $g(n)$ satisfies the second half of \eqref{e:g} by definition.
  By Lemma~\ref{l:xi}(i)
  \begin{equation}
    g(n)\le
    \min\Big\{k: K \binom n k^{-1} n^{1/2} \log n\le 2^{-n}m'(n)\Big\}
  \end{equation}
  Write $k=n/2-i$. Using formula \eqref{e:binomass} for the binomial
  coefficient we find that
  \begin{equation}
    g(n)\le \frac n2
    +\min\Big\{i: e^{2i^2/n}\le \frac {c m'(n)}{n\log n}\Big\}.
  \end{equation}
  Since $m'(n)\gg n\log n$ there exists $i_0$, $n^{2/3}\gg i_0\gg n^{1/2}$, for
  which the inequality in braces holds. Since $n/2-g(n)\ge i_0$, the first
  half of \eqref{e:g} is proved.
\end{proof}

Let us assume from now on that $g(n)$ satisfies \eqref{e:g}. As
announced earlier, the key ingredient of the proof of
Proposition~\ref{p:sum} is a precise estimate on the Laplace
transform of the hitting time of a single point. We now state and
prove this result. As we will see, this is where the function
$\xi_n$ comes in.

\begin{lemma}
  \label{l:laplace}
  Let $\boldsymbol 0$ be the vertex
  of the hypercube with all coordinates equal to $0$ and
  let $z_k\in \mathcal V_n$ be an arbitrary
  vertex of the hypercube such that $d(z_k,\boldsymbol 0)=k$. If
  $n\log n \ll m(n) \ll 2^n$, then for all $s>0$
  \begin{equation}
    \mathbb E_{z_k}\exp\Big(-\frac s {m(n)} H_n(\boldsymbol 0)\Big)=
    \Big[2^{-n}\frac {m(n)} s + \xi_n(k)\Big](1+o(1)).
  \end{equation}
\end{lemma}
\begin{proof}
  By Fourier methods for random walks on finite groups \cite{Dia88}, we
  have as in \cite{Mat89,BC06}
  \begin{equation}
    \label{e:furt}
    \mathbb E_{z_k} e^{-s/m H_n(\boldsymbol 0)}=
    \frac
    {\sum_{y\in \mathcal V_n}(-1)^{z_k\cdot y}\big[1-e^{-s/m }
    (1-\frac {2 d(y,\boldsymbol 0)}{n})\big]^{-1}}
    {\sum_{y\in \mathcal V_n}\big[1-e^{-s/m}
    (1-\frac {2 d(y,\boldsymbol 0)}{n})\big]^{-1}},
  \end{equation}
  where $x\cdot y = \sum_{i=1}^n x(i) y(i)$ is the standard scalar
  product in $\mathbb R^n$.

  Let us first consider the numerator of \eqref{e:furt}. Observe that there are
  $\tbinom k i \tbinom {n-k} j$ sites $y\in \mathcal V_n$ such that
  $d(0,y)=i+j$ and $z_k\cdot y = i$. Hence the numerator of \eqref{e:furt} is equal to
  \begin{equation}
    \sum_{i=0}^k \sum_{j=0}^{n-k}
    (-1)^i \binom k i \binom {n-k} j
    \frac 1 {1-e^{-s/m}(1-2n^{-1}(i+j))}.
  \end{equation}
  This expression can be simplified using the following lemma.
  \begin{lemma}
    For all $k,j\in \{0,1,\dots\}$ and all $s>0$
    \begin{equation}
      \label{e:betaint}
      \sum_{i=0}^k
      \frac{(-1)^i \binom k i}
      {1-e^{-s/m}(1-2n^{-1}(i+j))}
      =
      \frac {ne^{s/m}} 2 \cdot
      \frac{\Gamma (1+k)\Gamma(j+\frac n2(e^{s/m}-1))}
      {\Gamma(1+k+j+\frac n2 (e^{s/m}-1))}
      .
    \end{equation}
  \end{lemma}
  \begin{proof}
    Note that the second fraction on the right-hand side of
    \eqref{e:betaint} can be expressed using the Beta-integral,
    \begin{equation}
      \frac {ne^{s/m}} 2 \cdot
      \frac{\Gamma (1+k)\Gamma(j+\frac n2(e^{s/m}-1))}
      {\Gamma(1+k+j+\frac n2 (e^{s/m}-1))}
      =
      \frac {ne^{s/m}} 2
      \int_0^1
      (1-t)^k t^{i+\frac n2 (e^{s/m}-1)-1}\,\d t.
    \end{equation}
    Expanding $(1-t)^k$ according to the binomial theorem, and performing an easy
    integration then gives the left-hand side of \eqref{e:betaint}.
  \end{proof}

  Using the last lemma the numerator of \eqref{e:furt} can be rewritten as
  \begin{equation}
    \label{e:num2}
    \sum_{j=0}^{n-k}
    \binom {n-k} j
    \frac {ne^{s/m}} 2 \cdot
    \frac{\Gamma (1+k)\Gamma(j+\frac n2(e^{s/m}-1))}
    {\Gamma(1+k+j+\frac n2 (e^{s/m}-1))}.
  \end{equation}

  So far we obtained an exact expression which we now want to evaluate. To do so
  we will use the following known properties of $\Gamma $-functions
  (we refer to \cite{AS72} for the definition and the properties of the functions
  appearing below).
  \begin{equation}
    \Gamma '(x)=\Gamma (x) \psi _0(x),\qquad
    \Gamma''(x)=\Gamma (x)(\psi _0(x)^2+\psi _1(x)),\qquad
    \lim_{x\to 0}x\Gamma (x)=1,
  \end{equation}
  where $\psi _0$ is the digamma function and $\psi _1=\psi_0'$.
  The values of these functions for integer arguments can be written explicitly:
  \begin{equation}
    \psi _0(k)=\sum_{i=1}^{k-1} \frac 1 i - \gamma_E = \log k(1+o(1)), \qquad
    \psi _1(k)=\frac {\pi^2} 6- \sum_{i=1}^{k-1} \frac 1 {i^2}=O(1/k).
  \end{equation}
  where $\gamma_E$ is Euler's constant.

  We can now evaluate \eqref{e:num2}. Set $\varepsilon = \frac n2(e^{s/m}-1)$
  and observe that the bound $m\gg n\log n$ entails $\varepsilon \ll (\log n)^{-1}$.
  To treat the term $j=0$ in \eqref{e:num2} simply note that,
  since $\varepsilon \psi_0(1+k)=o(1)$ for all $k\le n$,
  \begin{equation}
    \label{e:blams}
    \frac {ne^{s/m}} 2 \cdot
    \frac{\Gamma (1+k)\Gamma(\varepsilon )}
    {\Gamma(1+k+\varepsilon )}
    =
    \frac {ne^{s/m}} 2 \cdot
    \frac {\varepsilon^{-1}(1+o(1))}
    {1+\varepsilon \psi _0(1+k)}=\frac m s
    (1+o(1)).
  \end{equation}
  A similar calculation for the remaining terms combined with
  \eqref{e:secondform} readily yields
  \begin{equation}
    \begin{split}
      \label{e:blaxi}
      \sum_{j=1}^{n-k}&
      \binom {n-k} j
      \frac {ne^{s/m}} 2 \cdot
      \frac{\Gamma (1+k)\Gamma(j+\varepsilon )}
      {\Gamma(1+k+j+\varepsilon )}
      \\&=
      \sum_{j=1}^{n-k}
      \frac n {2j} \binom {n-k} j {\binom{k+j} j}^{-1}
      \frac {(1+\varepsilon \psi _0(j))(1+o(1))}
      {1+\varepsilon \psi _0(1+k+j)}
      =2^n \xi_n(k)(1+o(1)).
    \end{split}
  \end{equation}

  For $k=0$ the numerator of \eqref{e:furt} coincides with the denominator.
  Equations \eqref{e:blams}, \eqref{e:blaxi} and Lemma~\ref{l:xi}(iv)
  then imply that the denominator behaves like
  \begin{equation}
    \label{e:bladem}
    \{ms^{-1} + 2^n \xi_n(0)\}(1+o(1))=2^n(1+o(1)).
  \end{equation}
  Finally, putting together \eqref{e:blams}, \eqref{e:blaxi} and \eqref{e:bladem}
  yields the claim of Lemma~\ref{l:laplace}.
\end{proof}

Lemma~\ref{l:laplace} now allows us to get information on the form
of the probability distribution function of $H_n(\boldsymbol 0)$.
Let us denote by $p_n(a,k)$ the probability
\begin{equation}
  p_n(a,k)=\mathbb P_{z_k}[H_n(\boldsymbol 0)< a m(n)].
\end{equation}

\begin{lemma}
  \label{l:twopoint}
  (i) There exists $C<\infty$ independent of $n$ and $k$ such that
  \begin{equation}
    \label{e:pnupperbound}
    p_n(a,k)\le Ce^a\big(2^{-n}m(n) +  \xi_n(k)\big).
  \end{equation}

  (ii)  For
  any $a\in [0,\infty)$, uniformly on compact subsets of this interval,
  \begin{equation}
    \lim_{n\to\infty}
    \max_{k\ge g(n)}
    \big |2^n m(n)^{-1}p_n(a,k)-a\big|=0.
  \end{equation}
\end{lemma}

\begin{proof}
  Assertion (i)  follows from Chebyshev
  inequality and Lemma~\ref{l:laplace}. To prove (ii) observe
  that
  for $k\ge g(n)$, by \eqref{e:g} and by Lemma~\ref{l:xi}(iii),
  $\xi_n(k)\le \xi_n(g(n))\ll 2^{-n}m(n)$. Therefore,
  $\mathbb E_{z_k}\exp(-s H_n(\boldsymbol 0)/m(n))=2^{-n}m(n)/s(1+o(1))$.
  Consider the sequence of measures $\mu_n$ given by
  \begin{equation}
    \mu_n([0,t])= 2^n m(n)^{-1} \mathbb P_{z_k}\big[H_n(\boldsymbol 0)/m(n)\in
      [0,t]\big].
  \end{equation}
  The Laplace transform of $\mu_n$ then satisfies
  $  \int_0^\infty e^{-s t}\mu_n(\d t)\xrightarrow{n\to\infty}  1/s$.
  Therefore, $\mu_n$ converges weakly to the Lebesgue measure (see
    \cite{Fel71},  Section XIII.1, Theorem 2a, p. 433) and thus
  \begin{equation}
    2^n m(n)^{-1}p_n(a,k)=\mu_n([0,a))\xrightarrow{n\to\infty} a.
  \end{equation}
  The uniformity on compact sets follows easily from the fact that
  the probabilities $p_n(a,k)$ are increasing in $a$.
\end{proof}

We finally use Lemma~\ref{l:twopoint} to get information on the
form of the probability distribution function of the hitting time
of finite subsets of points of $\mathcal V_n$. For
$y_1,\dots,y_\ell\in \mathcal V_n$ we define $\bar
H_n(y_1,\dots,y_\ell)$ by
\begin{equation}
  \bar H_n(y_1,\dots,y_\ell)=
  \begin{cases}
    H_n(y_\ell)&\text{if $H_n(y_i)<H_n(y_{i+1})$ for all
      $i\in \{1,\dots,\ell -1\}$},\\
    \infty&\text{otherwise}.
  \end{cases}
\end{equation}
That is $\bar H_n(y_1,\dots,y_\ell)$ is finite only if the $y$'s were
visited in the prescribed order. In this case it is equal to the time to
visit all $y$'s. Observe that it is always infinite if $y_i=y_j$ for some
$i\neq j$.

\begin{lemma}
  \label{l:coverbad}
  Let $x=y_0$, let $y_1,\dots,y_\ell$ be mutually distinct points in
  $\mathcal V_n$, and let $d(i)=d(y_{i-1},y_i)$. Then
  \begin{equation}
    \label{e:coverbad}
    \mathbb P_x\big[\bar H_n(y_1,\dots,y_\ell)<am(n)\big]
    \le
    C^\ell e^{\ell a} \prod_{i=1}^\ell \big(2^{-n}m(n)+ \xi_n(d(i))\big).
  \end{equation}
\end{lemma}

\begin{proof}
  Obviously, by the strong Markov property,
  \begin{equation}
    \mathbb P_x\big[\bar H_n(y_1,\dots,y_\ell)<am(n)\big] \le
    \prod_{i=1}^\ell
    \mathbb P_{y_{i-1}}[H_n(y_i)< am(n)],
  \end{equation}
  which, by Lemma~\ref{l:twopoint}(i), is bounded by the right-hand side of
  \eqref{e:coverbad}.
\end{proof}

\begin{lemma}
  \label{l:covergood}
  Let $x^n=y_0^n$, let $y_1^n,\dots,y_\ell^n\in \mathcal V_n$, and let
  $d_n(i)=d(y^n_i,y^n_{i-1})$. Suppose that $d_n(i)\ge g(n)$ for
  all $i\in\{1,\dots,\ell\}$ and all $n$.
  Then, uniformly over all $x^n$ and $y^n_i$,
  \begin{equation}
    \lim_{n\to\infty}
    2^{\ell n} m(n)^{-\ell}
    \mathbb P_{x^n}\big[\bar H_n(y_1^n,\dots,y_\ell^n)
      < a m(n)\big] =
    \frac {a^\ell}{\ell !}.
  \end{equation}
\end{lemma}

\begin{proof}
  The probability in question can be bounded from above using the strong Markov
  property,
  \begin{equation}
    \label{e:iiint}
    \mathbb P_{x^n}\big[\bar H_n(y_1^n,\dots,y_\ell^n)
      < a m(n)\big] \le
    \sum_{\substack{a_1,\dots, a_\ell\in \mathbb N/m(n)
        \\ a_1+\dots+a_\ell < a}}
    \prod_{i=1}^\ell
    \mathbb P_{y_{i-1}}[H_n(y_i)/m(n)= a_i].
  \end{equation}
  Since $d(y_{i-1},y_i)>g(n)$, it is easy to
  see from Lemma~\ref{l:twopoint}(ii) that the sum behaves like
  \begin{equation}
    2^{-n\ell} m(n)^{\ell}
    \idotsint\limits_{x_1+\dots+x_\ell < a}
    \d x_1\dots\d x_\ell
    (1+o(1))
    =
    2^{-n\ell} m(n)^{\ell}
    a^\ell / \ell !
    (1+o(1)).
  \end{equation}

  The above expression only provides an upper bound since it does not
  exclude the possibility that the random walk visits another $y_j$, $j>i$, on its way from
  $y_{i-1}$ to $y_i$.
  To construct a lower bound we should exclude such
  visits. Therefore, denoting by $\mathrm{UB}$ the upper
  bound~\eqref{e:iiint},
  \begin{equation}
    \begin{split}
      \mathbb P_{x^n}&\big[\bar H_n(y_1^n,\dots,y_\ell^n)
        < a m(n)\big]
      \\
      &\ge \mathrm{UB}-
      \sum_{i=1}^\ell
      \sum_{j=i+1}^\ell
      \sum_{\substack{a_1,\dots, a_\ell\in \mathbb N/m(n)
          \\ a_1+\dots+a_\ell < a}}
      \prod_{\substack{k=1\\k\neq i}}^\ell
      \mathbb P_{y_{k-1}}[H_n(y_k)/m(n)= a_k]
      \\&\quad\times
      \sum_{\substack{b\in \mathbb N/m(n)\\ b\in(0,a_i-a_{i-1})}}
      \mathbb P_{y_{i-1}}[H_n(y_j)/m(n)= b]
      \mathbb P_{y_{j}}[H_n(y_i)/m(n)= a_i-b].
    \end{split}
  \end{equation}
  The negative term on the right-hand side is smaller than
  \begin{equation}
    C2^{-n(\ell-1)} m(n)^{\ell-1}
    e^{(\ell-1)a}
    \sum_{i=1}^\ell
    \sum_{j=i+1}^\ell
    \mathbb P_{y_{i-1}}[H_n(y_j)< am(n)]
    \mathbb P_{y_{j}}[H_n(y_i)<am(n) ].
  \end{equation}
  Hence, if we can show that for all $j\neq i$
  \begin{equation}
    \label{e:smalll}
    \mathbb P_{y_{i-1}}[H_n(y_j)< am(n)]
    \mathbb P_{y_{j}}[H_n(y_i)<am(n) ] = o(2^{-n}m(n)),
  \end{equation}
  then the proof of Lemma~\ref{l:covergood} is finished.

  Let $k=d(y_{i-1},y_j)$ and $l=d(y_j,y_i)$.
  Since $d(y_{i-1},y_i)\ge g(n)$ we have also $k+l\ge g(n)$. Obviously
  $k\ge 1$, $l\ge 1$. By Lemma~\ref{l:twopoint}(i),
  \begin{equation}
    \begin{split}
      \mathbb P_{y_{i-1}}&[H_n(y_j)< am(n)]
      \mathbb P_{y_{j}}[H_n(y_i)<am(n) ]
      \\&\le
      C^2e^{2a}\big(2^{-n}m(n)+ \xi_n(k)\big).
      \big(2^{-n}m(n)+ \xi_n(l)\big)
      \\&\le
      C^2e^{2a}\big(2^{-2n}m(n)^{2}+ 2^{-n}m(n) 2\xi_n(1)+\xi_n(k)\xi_n(l)\big).
    \end{split}
  \end{equation}
  The first two summands are $o(2^{-n}m(n))$ (see Lemma~\ref{l:xi}(iii,iv)).
  If $k$ or $l$ is larger than $g(n)$, then the same is valid for the
  third one. As the last step of the proof we show that if
  $\max \{k,l\}< g (n)$ and $k+l\ge g(n)$, then for any $\varepsilon >0$ and
  $n$ large enough
  \begin{equation}
    \label{e:submult}
    \xi_n(k)\xi_n(l)\le \varepsilon 2^{-n} m(n).
  \end{equation}
  Let $z_{k+l}$ be as in Lemma~\ref{l:laplace}
  and let $z_k$ be any point such that $d(z_k,\boldsymbol 0)=k$ and
  $d(z_k,z_{k+l})=l$. Since on the way from $z_{k+l}$ to $\boldsymbol 0$,
  the random walk may pass through $z_k$ we have
  \begin{equation}
    \mathbb E_{z_{k+l}}[e^{-sH_n(\boldsymbol 0)}]
    \ge
    \mathbb E_{z_{k+l}}[e^{-sH_n(z_k)}]
    \mathbb E_{z_{k}}[e^{-sH_n(\boldsymbol 0)}].
  \end{equation}
  Lemma~\ref{l:laplace} then yields
  \begin{equation}
    \big\{2^{-n}\frac {m(n)} s + \xi_n(k+l)\big\}(1+o(1))
    \ge
    \big\{2^{-n}\frac {m(n)} s + \xi_n(k)\big\}
    \big\{2^{-n}\frac {m(n)} s + \xi_n(l)\big\}
    \ge \xi_n(k)\xi_n(l).
  \end{equation}
  Since $k+l\ge g(n)$ we can, in view of \eqref{e:g}, ignore the term
  $\xi_n(k+l)$ on the left-hand side. Taking $s$ sufficiently large then
  proves \eqref{e:submult}. This concludes the proof of the lemma.
\end{proof}

We are now ready to complete the proof of Proposition~\ref{p:sum}.

\subsection{Proof of Proposition~\ref{p:sum}}
\label{s:propproved}

  We shall establish that
  \begin{equation}
    \lim_{n\to\infty}
    \max_{x\in \mathcal V_n}
    \bigg|
    \psum_{y_1,\dots,y_i\in B_n\setminus x}
    \mathbb P_x\Big[\bigcap_{j=1}^i H_n(y_j)< a m(n)\Big]
    -a^i\bigg| =0.
  \end{equation}
  Observe that, with $y_0=x$,
  \begin{equation}
    \begin{split}
      \label{e:longsum}
      \psum_{y_1,\dots,y_i\in B_n\setminus x}
      \mathbb P_x\Big[&\bigcap_{j=1}^i H_n(y_j)< a m(n)\Big]
      =
      i!
      \psum_{y_1,\dots,y_i\in B_n\setminus x}
      \mathbb P_x\big[\bar H_n(y_1,\dots,y_i)< a m(n)\big]
      \\&=
      i! \sum_{d_1,\dots,d_i=1}^n
      \psum_{\substack{y_1,\dots,y_i\in B_n\\d(y_i,y_{i-1})=d_i}}
      \mathbb P_x\big[\bar H_n(y_1,\dots,y_i)< a m(n)\big].
    \end{split}
  \end{equation}
  Consider first the summation over distances larger than $g(n)$. Using
  Lemma~\ref{l:covergood} we get that (uniformly in the starting position $x$)
  \begin{equation}
    \begin{split}
      i! \sum_{d_1,\dots,d_i=g(n)}^n&
        \psum_{\substack{y_1,\dots,y_i\in B_n\\d(y_i,y_{i-1})=d_i}}
        \mathbb P_x\big[\bar H_n(y_1,\dots,y_i)< a m(n)\big]
        \\&=
        \sum_{d_1,\dots,d_i=g(n)}^n
        \psum_{\substack{y_1,\dots,y_i\in B_n\\d(y_i,y_{i-1})=d_i}}
        2^{-in} m(n)^i a^i(1+o(1))
        =a^i(1+o(1)).
      \end{split}
  \end{equation}
  For the second equality we used the fact that by \eqref{e:Bsize},
  \eqref{e:vcond} and the finiteness of $i$ there are
  $2^n m(n)^{-1}(1+o(1))$ choices for every $y_i$.

  To estimate the remaining contribution to \eqref{e:longsum},
  we first bound  the sum
  \begin{equation}
    \label{e:smalls}
    \sum_{d=1}^{g(n)-1}
    \sum_{y\in B_n:d(y,x)=d}
    C e^a  \big(2^{-n}m(n) + \xi_n(d)\big)
    \le
    \sum_{d=1}^{g(n)-1} v_n(d)\xi_n(d)
    +V_n(g(n)-1)  2^{-n}m(n).
  \end{equation}
  Both summands of in the last formula converge to $0$ which can be seen
  easily from  \eqref{e:vcond} and  \eqref{e:Bsize}. For $d>g(n)$,
  by Lemma~\ref{l:xi}(iii), $\xi_d(n)\ll 2^{-n}m(n)$.
  Therefore, using \eqref{e:Bsize}, for all $n$ large enough
  \begin{equation}
    \label{e:bigs}
    \sum_{d=1}^{n}
    \sum_{y\in B_n:d(y,x)=d} C e^a\big(2^{-n}m(n) + \xi_n(d)\big)
    \le 2Ce^a|B_n|2^{-n}m(n)\le 4 C e^a.
  \end{equation}
  According to  Lemma~\ref{l:coverbad},  the remaining
  part of \eqref{e:longsum} then  satisfies
  \begin{equation}
    \begin{split}
      i!& \sum_{\substack{d_1,\dots,d_i=1\\\exists d_i\le g(n)}}^n
      \psum_{\substack{y_1,\dots,y_i\in B_n\\d(y_i,y_{i-1})=d_i}}
      \mathbb P_x\big[\bar H_n(y_1,\dots,y_i)< a m(n)\big]
      \\&\le
      i!i
      \sum_{d_1=1}^{g(n)}
      \sum_{\substack{y\in B_n\\d(y,x)=d_1}}
      C a \big(2^{-n}m(n)+ \xi_n(d_1)\big)
      \bigg( \sum_{d=1}^{n}
        \sum_{\substack{y\in B_n\\d(y,x)=d}}
        C a\big(2^{-n}m(n)+ \xi_n(d)\big)
        \bigg)^{i-1},
      \end{split}
  \end{equation}
  which converges to $0$ by \eqref{e:smalls} and \eqref{e:bigs}. This
  finishes the proof of Proposition~\ref{p:sum}.
%%%% END OF ht2.tex %%%%

%%%% START OF ht3.tex %%%%
\section{Proof of Theorem~\ref{t:perc} and of Theorem~\ref{t:sampling}}
\label{s:equiv}

In this section we apply Theorem~\ref{t:gen} to derive the
asymptotic hitting distribution of randomly chosen sets in two
different settings: for random clouds (namely we prove
Theorem~\ref{t:perc}) and in the setting of drawing without
replacement (which is Theorem~\ref{t:sampling}).

\begin{proof}[Proof of  Theorem~\ref{t:perc}]
To prove Theorem~\ref{t:perc} we will naturally show that the
assumptions of
Theorem~\ref{t:gen} are satisfied for percolation clouds $A_n$ of
density $\bar m(n)^{-1}$, where $n\log n\ll \bar m(n)\ll 2^n(\log n)^{-1}$.

We first verify condition \eqref{e:Bsize}, i.e.~that
$P$-a.s.~$|A_n|=2^n \bar m(n)^{-1}(1+o(1))$.
By Chebyshev exponential inequality, for any $\delta >0$
with $\lambda >0$,
\begin{equation}
  \begin{split}
    \label{e:cheb}
     P[|A_n|\gtrless(1\pm \delta )2^nn^{-1}]
    &\le
    \exp\{\mp\lambda  (1\pm \delta )2^n\bar m(n)^{-1}\}
    (1+(e^{\pm \lambda }-1)\bar m(n)^{-1})^{2^n}
    \\&\le
    \exp\{\mp\lambda  (1\pm \delta )2^n\bar m(n)^{-1} +
      (e^{\pm \lambda }-1)2^n\bar m(n)^{-1}\}.
  \end{split}
\end{equation}
Taking $\lambda $ sufficiently small and using the fact that
$2^n\bar m(n)^{-1}\gg \log n$, we see that  the right-hand side of the last
equation is summable. Borel-Cantelli lemma then yields the result.

Let $f_n(k)=n(\log \bar m(n))^{-1}\bbone\{k=1\} + n\bbone\{k>1\}$. To prove
that the first part of \eqref{e:vcond} is satisfied for the percolation
cloud we show:
\begin{lemma}
  \label{l:cloudv}
  There exists $C$ large enough, such that for a.e.~realisation of $A_n$
  and for $n$ large enough
  \begin{equation}
    \label{e:rhs}
    v_n(k)\le C\Big[\binom n k \bar m(n)^{-1} +  f_n(k)\Big]\qquad\forall
    k\in\{1,\dots,g(n)\}.
  \end{equation}
\end{lemma}
\begin{proof}
  Let $F_n(k)$ denote the right-hand side of
  \eqref{e:rhs}. By definition of $v_n(k)$,
  \begin{equation}
    \begin{split}
      P[v_n(k)\ge F_n(k)]
      &\le
      \sum_{x\in \mathcal V_n}
      P[|\{y\in A_n:d(x,y)=k\}|\ge F_n(k)]
      \\&=
      2^n
      P[|\{y\in A_n:d(\boldsymbol 0,y)=k\}|\ge F_n(k)].
    \end{split}
  \end{equation}
  Using the same calculation as in \eqref{e:cheb} this is bounded from
  above by
  \begin{equation}
    2^n \exp\Big\{-C\lambda \Big[\binom n k \bar m(n)^{-1} + f_n(k)\Big]+
      (e^\lambda -1)\binom n k \bar m(n)^{-1}\Big\}.
  \end{equation}
  If we choose $\lambda = \log \bar m(n)$ and $C$ large enough for $k=1$, or
  $\lambda = \mathit{const}$ and $C$ large enough for $k>1$, then the
  right-hand side of the last equation decays at least as fast as $e^{-cn}$ for all
  $k$. Summing over $k$ and using Borel-Cantelli Lemma yields the desired
  result.
\end{proof}

Lemma~\ref{l:cloudv} implies that
\begin{equation}
  \sum_{k=1}^{g(n)-1} v_n(k) \xi_n(k)
  \le
  \sum_{k=1}^{g(n)-1}  C\Big[\binom n k \bar m(n)^{-1}+ f_n(k)\Big] \xi_n(k).
\end{equation}
Using Lemma~\ref{l:xi}(i),(iv) this can be bounded by
\begin{equation}
  \label{e:condcond}
  C\bigg\{
    \sum_{k=1}^{g(n)-1}
    \xi_n(k)\binom n k \bar m(n)^{-1}+
    \sum_{k=1}^3 f_n(k)\binom n k^{-1}
    +\sum_{k=4}^{g(n)-1} n  f_n(k) \binom n k^{-1}
    \bigg\}
\end{equation}
The last two terms in the last formula are bounded by
\begin{equation}
    C \sum_{k=1}^3 f_n(k) n^{-k}
    + C n^3 n^{-4}\xrightarrow{n\to\infty}0
\end{equation}
as can be seen easily from the definition of $f_n(k)$.
The first term in \eqref{e:condcond} equals
\begin{equation}
  \frac {n2^{-n}}{2\bar m(n)}
  \sum_{k=1}^{g(n)-1}
  \sum_{j=1}^{n-k}
  \binom n { j+k} \frac 1j
  =
  \frac {n}{2\bar m(n)}
  \sum_{j=1}^{n-1}
  \frac {2^{-n}}j
  \sum_{k=j+1}^{n\wedge (g(n)+j-1)}\binom n k\le
  \frac {C n \log n }{\bar m(n)}.
\end{equation}
This tends to $0$ by the assumptions on $\bar m(n)$.
Therefore $v_n(k)$ verifies the first part of \eqref{e:vcond} $P$-a.s.~.

To verify the second part observe first that if
$\bar m(n)\ge \delta 2^n n^{-1}$ for some $\delta >0$, then \eqref{e:gg}
holds for $g(n)=2$.  Therefore,
by Lemma~\ref{l:cloudv}, $V_n(g(n)-1)\le v_n(1) + 1 \le C \ll |A_n|$.
We can hence
further suppose that $\bar m(n)\ll 2^n n^{-1}$.
By moderate deviations argument, and since $n/2-g(n)\gg n^{1/2}$,
\begin{equation}
  |\{y:d(\boldsymbol 0,y)\le g(n)-1\}|\le 2^n e^{-c g(n)^2/n}.
\end{equation}
Since $2^n\bar m(n)^{-1}\gg n$ there is a function $f(n)$ such that
\begin{equation}
  2^n \bar m(n)^{-1}f(n)\gg n
  \qquad \text{and}\qquad
  1\gg f(n)\gg e^{-cg(n)^2/n}.
\end{equation}
As in \eqref{e:cheb}
\begin{equation}
  P\Big[V_n(g(n)-1)\ge
    \frac{2^{n}f(n)}{\bar m(n)}\Big]\le
  2^n\exp\Big\{\frac{2^n}{\bar m(n)}\big[
    -\lambda f(n)+
    (e^\lambda -1)e^{-cg(n)^2/n}
    \big]
    \Big\}.
\end{equation}
For our choice of $f$ this is summable. Therefore a.s.~for $n$ large enough
$V_n(g(n)-1)\le 2^n \bar m(n)^{-1} f(n)\ll |A_n|$. This verifies the second part of
\eqref{e:vcond}.

We have verified that with $P$ probability one the sequence of percolation
clouds $A_n$ satisfies all the assumptions of Theorem~\ref{t:gen}.
This proves Theorem~\ref{t:perc}.
\end{proof}

\begin{proof}[Proof of  Theorem~\ref{t:sampling}]
  For any $m(n)$ satisfying the conditions of Theorem~\ref{t:sampling} it
  is possible to choose $\bar m(n)$ satisfying the conditions of
  Theorem~\ref{t:perc} such that
  \begin{equation}
    \bar m(n)^{-1} \ge (1+\varepsilon )m(n).
  \end{equation}
  We now consider a sequence of percolation clouds $A_n$ with density
  $\bar m(n)^{-1}$ defined on the same probability space
  $(\Omega',F',P')$ as $A'_n$. Since $P'$-a.s.~for all $n$ large enough
  $|A_n|>2^n m(n)^{-1}=|A'_n|$, we can couple $A'_n$ and $A_n$ in the way
  that for all $n$ large $A'_n\subset A_n$. Moreover, $A_n$ satisfies the
  conditions \eqref{e:volumes}--\eqref{e:vcond} of Theorem~\ref{t:gen}. To
  finish the proof observe that if $A_n$ satisfies these conditions, then
  any subset of $A_n$ satisfies them too.
\end{proof}
%%%% END OF ht3.tex %%%%

%\input htapp

%%%% START OF htbiblio.tex %%%%
\subsection*{Acknowledgements} Both authors thank
the Chair of Stochastic Modelling
of the \'Ecole Polytechnique F\'ed\'erale of Lausanne
for financial support. Ji\v r\'\i~\v Cern\'y
thank the Centre de Physique Th\'eorique of Marseille for hospitality.

\bibliographystyle{amsalpha}
\bibliography{ht0}

\def\cprime{$'$}
\providecommand{\bysame}{\leavevmode\hbox to3em{\hrulefill}\thinspace}
\providecommand{\MR}{\relax\ifhmode\unskip\space\fi MR }
% \MRhref is called by the amsart/book/proc definition of \MR.
\providecommand{\MRhref}[2]{%
  \href{http://www.ams.org/mathscinet-getitem?mr=#1}{#2}
}
\providecommand{\href}[2]{#2}
\begin{thebibliography}{BEGK02}

\bibitem[AS72]{AS72}
Milton Abramowitz and Irene~A. Stegun, \emph{Handbook of mathematical functions
  with formulas, graphs, and mathematical tables.}, National Bureau of
  Standards Applied Mathematics Series, 55, fifth edition, John Wiley \& Sons
  Inc., Washington, D.C., 1972.

\bibitem[BB02]{BB02b}
E.~Bertin and J.-P. Bouchaud, \emph{Dynamical ultrametricity in the critical
  trap model}, J. Phys. A: Math. Gen. \textbf{35} (2002), 3039.

\bibitem[BBG03a]{BBG03}
G{\'e}rard Be{n A}rous, Anton Bovier, and V{\'e}ronique Gayrard, \emph{Glauber
  dynamics of the random energy model. {I}. {M}etastable motion on the extreme
  states}, Comm. Math. Phys. \textbf{235} (2003), no.~3, 379--425.

\bibitem[BBG03b]{BBG03b}
G{\'e}rard {Ben A}rous, Anton Bovier, and V{\'e}ronique Gayrard, \emph{Glauber
  dynamics of the random energy model. {II}. {A}ging below the critical
  temperature}, Comm. Math. Phys. \textbf{236} (2003), no.~1, 1--54.

\bibitem[B{\v C}06a]{BC06}
G{\'e}rard Be{n A}rous and Ji{\v r\'\i} {\v C}ern{\'y}, \emph{The arcsine law
  as a universal aging scheme for trap models}, to appear in Comm. Pure Appl.
  Math. (2006).

\bibitem[B{\v C}06b]{BC06c}
G{\'e}rard {Be{n A}rous} and Ji{\v r\'\i} {\v C}ern{\'y}, \emph{Dynamics of
  trap models}, {\'E}cole d'{\'et\'e} de physique des Houches, Session LXXXIII,
  Mathematical Statistical Physics, Elsevier, 2006, pp.~331--394.

\bibitem[BEGK01]{BEGK1}
Anton Bovier, Michael Eckhoff, V{\'e}ronique Gayrard, and Marcus Klein,
  \emph{Metastability in stochastic dynamics of disordered mean field models},
  Prob. Theor. Rel. Fields (2001), no.~119, 99--161.

\bibitem[BEGK02]{BEGK2}
Anton Bovier, Michael Eckhoff, V{\'e}ronique Gayrard, and Markus Klein,
  \emph{Metastability and low lying spectra in reversible {M}arkov chains},
  Comm. Math. Phys. \textbf{228} (2002), 219--255.

\bibitem[BG06]{BG06}
G{\'e}rard Be{n A}rous and V{\'e}ronique Gayrard, \emph{Elementary potential
  theory on the hypercube}, preprint {{\tt math.PR/0611178}}, 2006.

\bibitem[Dia88]{Dia88}
Persi Diaconis, \emph{Group representations in probability and statistics},
  Institute of Mathematical Statistics Lecture Notes---Monograph Series, 11,
  Institute of Mathematical Statistics, Hayward, CA, 1988.

\bibitem[Fel71]{Fel71}
William Feller, \emph{An introduction to probability theory and its
  applications. {V}ol. {II}.}, Second edition, John Wiley \& Sons Inc., New
  York, 1971.

\bibitem[Kem61]{Kem61}
J.H.B. Kemperman, \emph{The passage problem for a stationary markov chain},
  Satistical Reseach Monographs, Vol. I, University of Chicago Press, Chicago,
  1961.

\bibitem[Mat88]{Mat88}
Peter Matthews, \emph{Covering problems for {M}arkov chains}, Ann. Probab.
  \textbf{16} (1988), no.~3, 1215--1228.

\bibitem[Mat89]{Mat89}
Peter Matthews{}, \emph{Some sample path properties of a random walk on the
  cube}, J. Theoret. Probab. \textbf{2} (1989), no.~1, 129--146.

\end{thebibliography}
%%%% END OF htbiblio.tex %%%%

\end{document}